\documentclass{article}
\usepackage[utf8]{inputenc}
\usepackage{natbib}
\usepackage{booktabs}
\usepackage{amssymb,amsmath,amsthm}
\usepackage{dsfont} 
\usepackage{graphicx}
\usepackage[textwidth=3cm]{todonotes}

\graphicspath{{img/}}

\usepackage{hyperref}
\hypersetup{
    colorlinks=true,
    linkcolor=blue,
    filecolor=magenta,      
    urlcolor=cyan,
    pdfpagemode=FullScreen,
    }
    
\title{Cluster Size Matters: A Comparative Study of Notip and pARI for Post Hoc Inference in fMRI}
\author{Nils Peyrouset,$^{1}$ Pierre Neuvial$^{1\ast}$, Bertrand Thirion$^{2}$\\
{\small $^{1}$Institut de Mathématiques de Toulouse; Université de Toulouse; CNRS;}\\
{\small UPS, F-31062 Toulouse Cedex 9, France}\\
{\small $^{2}$INRIA, Université Paris-Saclay, CEA}\\
{\small $^\ast$Correspondence:  pierre.neuvial@math.univ-toulouse.fr}
}
\date{}

\begin{document}
\maketitle

\begin{abstract}
All Resolutions Inference (ARI) is a post hoc inference method for functional Magnetic Resonance Imaging (fMRI) data analysis that provides valid lower bounds on the proportion of truly active voxels within any, possibly data-driven, cluster.
As such, it addresses the paradox of spatial specificity encountered with more classical cluster-extent thresholding methods. 
It allows the cluster-forming threshold to be increased in order to locate the signal with greater spatial precision without overfitting, also known as the drill-down approach.
Notip and pARI are two recent permutation-based extensions of ARI designed to increase statistical power by accounting for the strong dependence structure typical of fMRI data. 

A recent comparison between these papers based on large voxel clusters concluded that pARI outperforms Notip.
We revisit this conclusion by conducting a systematic comparison of the two.
Our reanalysis of the same fMRI data sets from the Neurovault database  demonstrates the existence of complementary performance regimes: while pARI indeed achieves higher sensitivity for large clusters, Notip provides more informative and robust results for smaller clusters.
In particular, while Notip supports informative “drill-down” exploration into subregions of activation, pARI often yields non-informative bounds in such cases, and can even underperform the baseline ARI method.
\end{abstract}

\section{Introduction}

A classical approach to statistical inference for functional Magnetic Resonance Imaging (fMRI) data is cluster-extent-based thresholding. This method aims to identify clusters of adjacent voxels containing at least one active voxel \citep{nichols2003controlling}.
This approach suffers from two known limitations.
First, larger clusters provide less information than smaller ones, a phenomenon known as the spatial specificity paradox \citep{Woo2014}.
Second, when clusters are zoomed-in (or ``drilled-down'') by choosing a more stringent threshold, a form of double dipping occurs, resulting in the loss of statistical control \citep{kriegeskorte2009circular}.

Post hoc inference aims to address these limitations by providing statistical guarantees on the number or proportion of active voxels in arbitrary, possibly user-defined clusters \citep{goeman2011multiple}.
The first application of post hoc inference to fMRI data is the All Resolutions Inference (ARI) \citep{rosenblatt2018ari}. 
ARI relies on the Simes inequality \citep{simes1986improved}, which can be conservative due to the strong positive dependence typically found in fMRI data \cite{andreella2023permutation}.
Recently, several improvements to ARI have been proposed specifically for fMRI data analysis, leading to the Notip \citep{blain2022notip} and pARI \citep{andreella2023permutation} approaches.

All of these methods are based on the idea of calibration~\citep{blanchard2020post}, which uses permutation or sign-flipping to learn and adapt to the dependency structure of the data set at hand.
From a theoretical perspective, these methods differ only in the choice of the so-called \emph{template} (set of thresholds), which implicitly determines the relative weight given to smaller or larger sets of voxels.S
A recent comparison between Notip and pARI based on the data sets originally analyzed in \citep{blain2022notip} has shown that pARI outperforms Notip for large clusters \citep{andreella2024selective}.

In this paper, based on the same data sets, we demonstrate the existence of regimes in which one method outperforms the other one, a conclusion in line with the idea of ``No Free Lunch''.
The increased sensitivity of pARI for larger clusters (already noted by \cite{andreella2024selective}) comes at the price of decreased sensitivity for smaller clusters. 
In practice, for clusters of several hundreds of voxels, pARI is generally outperformed by Notip, but also by the baseline method ARI.
This implies that contrary to Notip and ARI, pARI does not provide informative results when drilling down into the clusters with the highest signal values.

In the remainder of the papier, we provide a self-contained description of the compared methods (Section~\ref{sec:methods}) and report extensive results on 37 fMRI data sets (Section~\ref{sec:results}), and provide a short discussion of the consequences of these results (Section~\ref{sec:discussion}).

\section{Methods}
\label{sec:methods}
\subsection{Post hoc inference for true discovery proportions}
\label{sec:goal-posthoc}

For each of $m$ voxels, we test the null hypothesis that voxel $i$ is not active under the condition of interest. 
The set $\mathcal{H}$ of all tested hypotheses is then identified to the set of all voxels, i.e. $\mathcal{H} = [m]$, where $[n] = \{1, \dots, n\}$ for any integer $n$. 
We denote by $\mathcal{H}_0 \subset \mathcal{H}$ the (unknown) subset of true null hypotheses. 
Let $m_0 = |\mathcal{H}_0 |$ be the (unknown) number of true null hypotheses and $\pi_0 = m_0 / m$ be the corresponding proportion.
For an arbitrary selection of voxels $S \subset \mathcal{H}$,  $|S \cap \mathcal{H}_0|$ is  the number of false positives within $S$, that is, the number of voxels that are selected whereas their corresponding null hypothesis is true (i.e., inactive voxels). 
The corresponding True Discovery Proportion (TDP) is then defined as $\mathrm{TDP}(S) = 1 - |S \cap \mathcal{H}_0|/|S|$.

Post hoc inference \citep{goeman2011multiple} aims at building an $(1-\alpha)$ a \textbf{TDP lower bound}, that is, a function $\overline{\mathrm{TDP}}_\alpha$ such that
\begin{align}
\label{eq:posthoc-TDP}
\mathbb{P} \bigg( 
\forall S \subset \mathcal{H},\:\:\: 
\mathrm{TDP}(S) \geq \overline{\mathrm{TDP}}_\alpha(S)
\bigg) \geq  1-\alpha .
\end{align}
That is, with probability greater than $1-\alpha$, the proportion of true discoveries of any subset $S$ is at least  $\overline{\mathrm{TDP}}_\alpha$. 
We emphasize that the "$for \, all \, S$" in \eqref{eq:posthoc-TDP} is \emph{inside} the probability: this implies that a TDP lower bound is valid for any number of possibly data-driven sets $S$. 
In the context of fMRI studies, such a TDP lower bound is applicable to all voxel clusters obtained by thresholding a statistical map.
Moreover, multiple cluster-forming thresholds may be chosen, possibly based on the results of the data analysis, without compromising the statistical validity of the TDP lower bound.
Therefore, as argued by \cite{rosenblatt2018ari}, post hoc methods address the problem of double dipping in fMRI data analysis, allowing users to  ``'drill down' from the cluster level to sub-regions, and even to individual voxels, in order to pinpoint the origin of the activation''.

\paragraph{Joint Error Rate Control.}
\cite{blanchard2020post} have shown that post hoc bounds may be sytematically derived from the control of a statistical risk called the Joint Error Rate, by a simple interpolation principle. 
We consider a  vector of $p$-values associated with each voxel: $\mathbf{p} = (p_1, \dots, p_m)$.
For a positive integer $K$, let $\mathbf{t} = (t_k)_{k \in [K]}$ be a non-decreasing vector of thresholds in $(0,1)$ aka template.  
The Joint Error Rate of the family $\mathbf{t}$ is defined by
\begin{align}
\mathrm{JER}(\mathbf{t}) 
= 
\mathbb{P} \big(\exists k \in \mathcal{H}_0 \cap [K], p_{(k:\mathcal{H}_0)} < t_k  \big),
\label{eq:JER}
\end{align}
where for $A \subset \mathcal{H}$ we denote by $p_{(k:A)}$ the $k$-th smallest $p$-value among $(p_i)_{i \in A}$.
\cite{blanchard2020post} have shown that if 
$\mathrm{JER}(\mathbf{t}) \leq \alpha$, then a TDP lower bound \eqref{eq:posthoc-TDP} is given by the function $\overline{\textrm{TDP}}^{\mathbf{t}}$, defined for $S \subset \mathcal{H}$ by 
\begin{align}
\overline{\textrm{TDP}}^{\mathbf{t}} (S)  
=
|S|^{-1} 
\left(
\max_{k \in [K]} 1 - k + \sum_{i \in S} \mathds{1} \{p_i < t_k\}
\right).
\label{eq:interpolation-bound}
\end{align}
After an initial sorting of the $p$-values, the bound $\overline{\mathrm{TDP}}^{\mathbf{t}}(S)$ can be computed in linear time ($O(|S|)$) from  \eqref{eq:interpolation-bound}, as shown by  \cite{enjalbert-courrech2022powerful}\footnote{A generic implementation applicable to any JER controlling family $\mathbf{t}$ is provided in the R package \texttt{sanssouci} and in the Python package \texttt{sanssouci.python}.}.
This framework for deriving post hoc bounds is particularly convenient in practice.
Indeed, since the \emph{computational} problem of the efficient evaluation of the post hoc bound is solved once and for all, the only remaining challenge to obtain a post hoc bound is the \emph{statistical} problem of finding a JER controlling family $\mathbf{t}$.

\subsection{Building JER controlling families}
\paragraph{The Simes family.}
The simplest example of JER controlling family is the Simes family, defined by $t_k=\alpha k /m$ for $k \in [m]$.
The \cite{simes1986improved} inequality states that for independent or positively associated test statistics \citep{sarkar2008simes}, in the sense of the PRDS property introduced by \cite{benjamini2001control}, we have:
\begin{align}
\mathbb{P} \big(\exists k \in \mathcal{H}_0, p_{(k:\mathcal{H}_0)} < \alpha k / m_0  \big) \leq
\alpha.
\label{eq:simes}
\end{align}
As the left hand side of \eqref{eq:simes} is exactly the JER of the Simes family, the Simes inequality trivially implies that the Simes family controls JER at level $m_0 \alpha/m = \pi_0 \alpha$, and a fortori at level $\alpha$.
It has been shown in \cite{blanchard2020post} that the post hoc bound obtained by interpolation recovers the Simes bound obtained by \cite{goeman2011multiple} by combining closed testing \citep{marcus1976closed} with a dedicated computational shortcut. 
%
The All Resolutions Inference (ARI) method \citep{rosenblatt2018ari} is an improved version of the above Simes-based method, where the thresholds $t_k= \alpha k /m $ are replaced with $t_k= \alpha k /h(\alpha)$, where $h(\alpha) \leq m$ is the Hommel value introduced in \cite{hommel1988stagewise}, which satisfies $m_0 \leq h(\alpha)$ with probabilty $1-\alpha$ \citep{goeman2019simultaneous}. 

The Simes inequality and the ARI method, which is based on it, are usually conservative in high-dimensional cases with dependent test statistics ~\citep{blanchard2020post,enjalbert-courrech2022powerful}. 
Such situations are a common use case in neuroimaging or genomic data (see e.g. \cite{Hayasaka2003-tr}).
This translates into the conservativeness of the associated post hoc bound: in such scenarios, the coverage of the post hoc bound \eqref{eq:posthoc-TDP} can be substantially larger than $1 - \alpha$.

\paragraph{JER Calibration.}
To address this conservativeness, a natural idea is to seek for other JER controlling families, whose JER is closer to the \emph{risk budget} $\alpha$. 
Given a threshold family $\mathbf{t}$ the JER \eqref{eq:JER} only depends on the \emph{joint} distribution of the null $p$-values, which is generally unknown.
To address this issue, \cite{blanchard2020post} have introduced a generic approach to approximating the JER. This approach involves sampling from the joint distribution using randomization-based methods.
In particular, the work of \cite{blanchard2020post} covers the classical cases of group label permutations for two-group testing and sign flipping for one-group testing.
More general linear models are covered in \cite{davenport2025fdp}.
Starting from a set of candidate families $\mathbf{t}$ called a template, JER calibration methods select the family  $\mathbf{t}^{\star}$ whose JER is the largest among those below the target risk/budget $\alpha$.
For a graphical illustration of the JER calibration principle, see \cite{blanchard2021agnostic,blain2022notip,andreella2024selective}.

\subsection{Existing JER calibration methods}

Following \cite{andreella2024selective}, we focus on the two most recent post hoc inference methods for the mass-univariate analysis of neuroimaging data \citep{blain2022notip,andreella2023permutation}. 
Both of them are based on JER calibration, and they only differ by the choice of candidate families (a.k.a. template).

\cite{andreella2023permutation} have introduced the permutation ARI (pARI) method. 
As it names suggests, it is inspired by the ARI method: the recommended candidate threshold families for neuroimaging data are of the form $t_k^{\delta}(\lambda) = \frac{\lambda(k-\delta)}{(m-\delta)} \mathds{1}_{\{k > \delta\}}$, where $\delta$ is an integer hyperparameter which has to be specified before data analysis to avoid circularity issues.
This hyperparameter indirectly controls where the method concentrates its power, via the minimal size of a region where non trivial inference can be made.
The choice $\delta=0$ recovers the Calibrated Simes method introduced by \cite{blanchard2020post}, whose numerical perfomance had already been studied by \cite{enjalbert-courrech2022powerful} for genomic applications. 
\cite{andreella2023permutation} recommend ``fixing $\delta = 1$ if the practitioner is interested in computing the lower bound for the TDP in small clusters, while $\delta > 1$ if the attention is focused on large cluster''.
In their application to fMRI data, \cite{andreella2023permutation} chose $\delta=1$ for their analysis of Auditory Data, and $\delta=27$ for their analysis of Rhyme Data. 
According to the follow-up paper \cite{andreella2024selective}, the choice $\delta=27$ is recommended for the analysis of fMRI data.

\cite{blain2022notip} have introduced the Notip method, where the main innovation is that the candidate threshold families are \emph{data-driven} instead of considering a pre-specified parametric part. 
In practice, Notip performs a first round of permutation on the data set at hand, and uses the successive empirical quantiles of the obtained null statistics as threshold families. 
The size $K$ of the threshold families is set to $2\%$ of the total number of voxels, that is, $K=1000$ when $m=50000$ voxels\footnote{Assuming that the proportion of active voxels in a typical fMRI data set is typically small (say, below $5\%$) and considering that TDP bounds below 1/2 are not informative, \cite{blain2022notip} have shown that one can focus on the $2.5\%$ percent of largest $p$-values, rounded to $2\%$ for simplicity.}.

\section{Results}
\label{sec:results}
Following \cite{andreella2024selective}, we consider both methods using the parameter values recommended by their authors, i.e. $\delta=27$ for pARI and $k_{\max}=1000$ for Notip. We start by studying one particular contrast, and then give general results on a set of 37 contrasts. 

\subsection{Focus on one contrast}
\label{sec:cluster-eval-one-contrast}

Here, we focus on the ``Look  negative  cue  vs  Look  negative  rating'' contrast, taken from the Neurovault collection\footnote{The corresponding data are available from  \url{http://neurovault.org/collections/1952}.}.
This dataset was already studied in \cite[Figures 5 and 7]{blain2022notip} to illustrate the face validity of the Notip method.
It was also used in \cite[Tables 2 to 5]{blain2022notip} to compare Notip to the baseline method ARI and to pARI with $\delta=0$ (refered to as ``calibrated Simes`` in \cite{blain2022notip}.  
The FDP bounds obtained by these methods are compared for each cluster obtained by a cluster-defining threshold of $z \in \{2.5, 3, 3.5\}$.
This comparison shows that Notip outperforms the other methods available at that time (therefore not including pARI with $\delta=27$), for all values of $z$.

\begin{figure}[!ht]
    \centering
    \includegraphics[width=\textwidth]{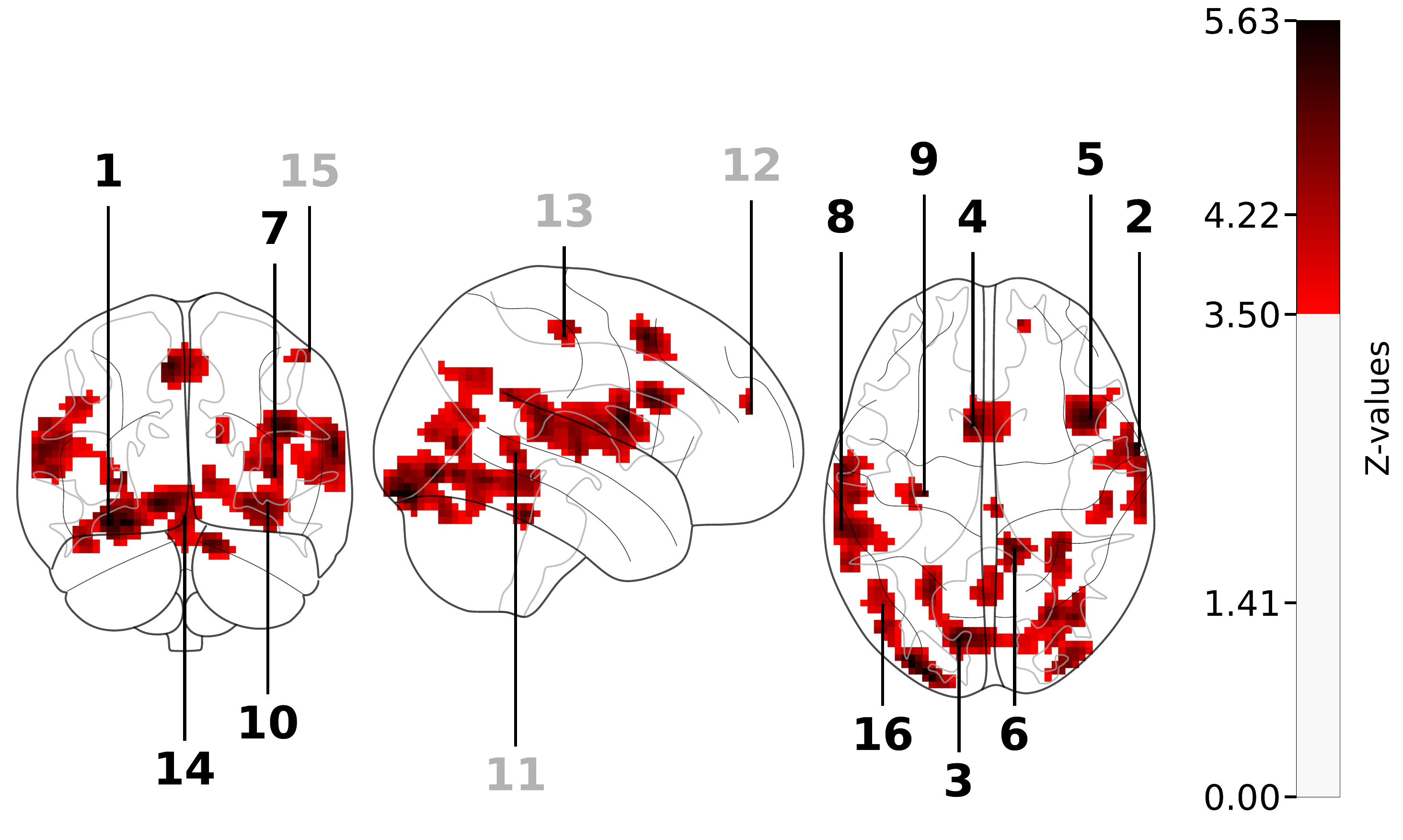}
    
\begin{tabular}{lrrrrrlll}
\toprule
\multicolumn{6}{c}{ } & \multicolumn{3}{c}{TDP lower bound} \\
\cmidrule(l{3pt}r{3pt}){7-9}
ID & X & Y & Z & Peak Stat & Size (mm3) & ARI & Notip & pARI\\
\midrule
1 & -33 & -94 & -17 & 5.63 & 3213 & 0.38 & \textbf{0.55} & 0.48\\
2 & 66 & 2 & 16 & 5.47 & 7425 & 0.38 & \textbf{0.77} & \textbf{0.77}\\
3 & -12 & -82 & -8 & 5.40 & 8397 & 0.46 & 0.79 & \textbf{0.8}\\
4 & -6 & 11 & 52 & 5.30 & 3321 & 0.23 & \textbf{0.5} & 0.49\\
5 & 45 & 14 & 25 & 5.27 & 2835 & 0.38 & \textbf{0.52} & 0.46\\
6 & 12 & -43 & -26 & 5.08 & 1107 & 0.15 & \textbf{0.2} & 0\\
7 & 39 & -73 & 4 & 5.00 & 2862 & 0.08 & \textbf{0.43} & 0.42\\
8 & -63 & -34 & 16 & 4.95 & 9585 & 0.46 & \textbf{0.82} & \textbf{0.82}\\
9 & -27 & -19 & 4 & 4.85 & 837 & \textbf{0.06} & \textbf{0.06} & 0\\
10 & 36 & -94 & -8 & 4.75 & 2160 & 0.25 & \textbf{0.42} & 0.3\\
14 & 0 & -64 & -14 & 4.43 & 1755 & 0 & \textbf{0.25} & 0.14\\
16 & -45 & -67 & 34 & 4.32 & 1890 & 0 & \textbf{0.21} & 0.13\\
\bottomrule
\end{tabular}

    \caption{Clusters identified with threshold $z = 3.5$ for the ``Look  negative  cue''  vs  ``Look  negative  rating'' data set: glass brain plot (top) and comparison between TDP lower bounds (bottom) For each cluster, the values in bold indicate the best result. Only clusters for which signal is detected by at least one method are reported.}
    \label{fig:TDP-cluster-table_task36_z=3.5}
\end{figure}

This comparison was complemented by \cite[Figure 4 and Table 1]{andreella2024selective}, where the pARI method with $\delta=27$ was added for $z=3$.
They observed that pARI (with $\delta=27$) outperformed Notip in this case. 
We were able to reproduce this observation (see Table \ref{tab:TDP-cluster-table_task36_z=3}). 
In order to complement this study, we have considered other choices for $z$. 
The results for $z=3.5$ are reported in Figure \ref{fig:TDP-cluster-table_task36_z=3.5}, where the clusters are represented on a glass brain plot.
The comparison results are more contrasted than those reported in \cite{andreella2024selective} for $z=3$: while Notip and pARI yield comparable TDP bounds for the largest clusters (clusters 2, 3 and 8), Notip performs better than pARI for smaller clusters. 
This behavior is somewhat expected: as noted by \cite{andreella2024selective}, ``$k_{\max}$
focuses power of Notip away from very large clusters, while $\delta>0$ focuses power of pARI away from small ones''.

\begin{figure}[!h]
    \centering
    \includegraphics[width=0.95\textwidth]{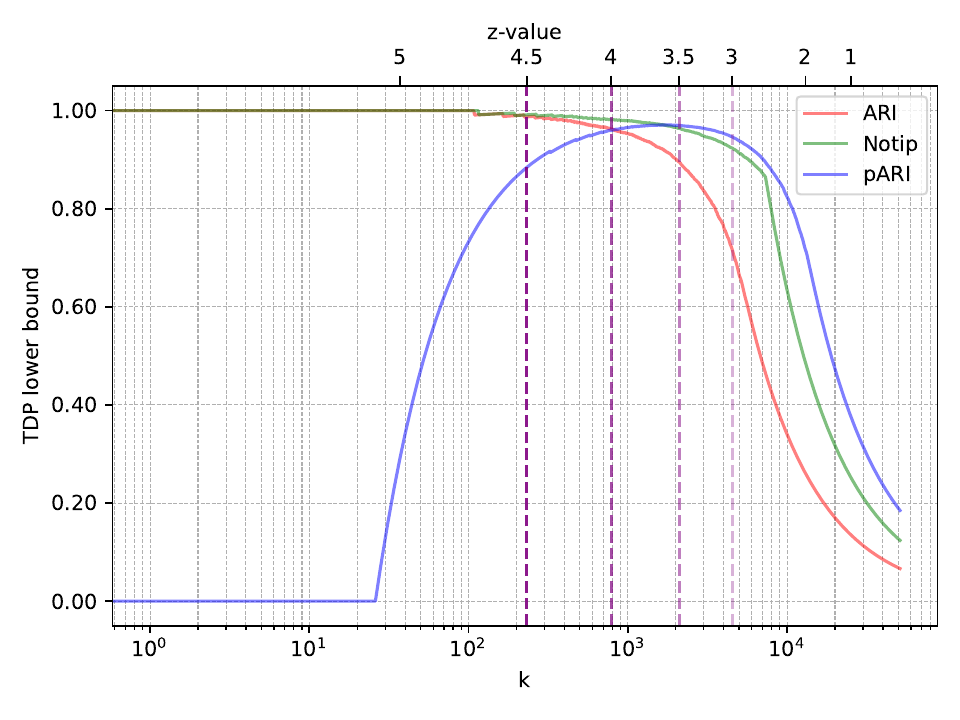}
    \caption{Confidence curve on the TDP for the ``Look  negative  cue  vs Look  negative  rating''  contrast: for each $k \in [m]$, we plot the TDP lower bound $\overline{\mathrm{TDP}}(S_k)$, where $S_k = \{i \in \mathcal{H}, |Z_i| \geq Z_{(k)}\}$ is the set of voxels with the $k$ largest $Z$ scores.}
    \label{fig:conf-curve}
\end{figure}

A more complete picture is brought by Figure~\ref{fig:conf-curve}, where the TDP lower bounds associated with all possible choices of $z$ (or equivalently, all possible $p$-value level sets) are displayed for each method.
For each value of $k$, we plot for each method the TDP lower bound $\overline{\mathrm{TDP}}(S_k)$ obtained for the set 
$S_k = \{i \in \mathcal{H}, |Z_i| \geq Z_{(k)}\}$
of voxels corresponding to the $k$ largest $Z$ scores.
In particular, the values of $k$ corresponding to $Z_{(k)} \in \{3, 3.5, 4, 4.5\}$ are highlighted by dotted vertical lines.

Since all of the compared method are valid post hoc bounds, the performance of a method can be quantified by the TDP bound, with higher values corresponding to a better bound. 
First, the Notip and pARI curves cross each other: this reflects the fact that no method is uniformly more powerful than the other one.
This point illustrates the absence of "free lunch" predicted by the theory outlined in Section \ref{sec:methods}: both methods optimize the same objective function, targeting a JER of $\alpha$ by estimating the joint null $p$-value distribution using permutations.
However, they have different constraints, which are encoded by the choice of a template.

As illustrated in Figure~\ref{fig:TDP-cluster-table_task36_z=3.5} for a specific value of $z$, the performance of Notip and pARI is comparable in the region $3 \leq z \leq 4$, with Notip better for larger values of $z$ (i.e., smaller values of $k$) and pARI better than Notip for smaller values of $z$ (i.e., larger values of $k$).
For smaller sets, the performance of pARI drops massively. 
This is expected for very small sets: by construction, pARI cannot detect any signal in sets of size less than $\delta=27$, corresponding here to $729$ mm$^3$.
However, pARI performs worse than the baseline ARI method for sets smaller than $800$ voxels.
This is alarming since the ARI method is known to be conservative for fMRI data \citep{blain2022notip}.
In fact, this conservativeness was the main motivation of the pARI method \citep{andreella2023permutation}.

%

A major feature of post hoc methods is their ability to ``further 'drill down' from the cluster level to sub-regions, and even to individual voxels, in order to pinpoint the origin of the activation'' \citep{rosenblatt2018ari}. 
In theory, all of the methods discussed here have this capacity since they provide TDP lower bounds that are valid for all possible sets of voxels simultaneously.
However, statistical validity (i.e., JER control) does not necessarily imply statistical power. 
In particular, \emph{the TDP lower bounds obtained by pARI (with $\delta=27$) for small sets of voxels are  non-informative}.
 Table~\ref{tab:TDP-cluster-table_task36_z=4} illustrates this point numerically. It provides the TDP lower bounds for the same dataset when drilling down to $z>4$:
pARI provides trivial (i.e. null) lower bounds for most clusters, and is outperformed by ARI (and a fortiori by Notip) even for the largest clusters of more than $3,000$ mm$^3$, corresponding to more than $100$ voxels. 
In practice, the pARI method cannot drill down to $z>4$ in this example.
In contrast, the Notip method is uniformly more powerful than the baseline ARI method and enables informative drilling down.
In the next section, we demonstrate that these observations are general, and not specific to this particular data set.

\subsection{fMRI datasets from the Neurovault database}

To consolidate the above results, we conducted experiments on a large fMRI data set: collection 1952 \citep{varoquaux2018atlases} of the Neurovault database (http://neurovault.org/collections/1952).
This dataset is an aggregation of 20 different fMRI studies and consists of statistical maps obtained at the individual level for a large set of contrasts.
We focused on 37 fMRI contrasts: the  ``Look  negative  cue  vs  Look  negative  rating'' contrast studied above, and the 36 contrasts introduced in \cite{blain2022notip} and further studied in \cite{andreella2024selective}.

We perform the same analysis for each contrast as in Section~\ref{sec:cluster-eval-one-contrast}.
For each threshold value of the $Z$ statistic and each contrast, we obtain a list of clusters and compute a TDP bound for each compared method. 
The results corresponding to Table \ref{tab:TDP-cluster-table_task36_z=3} and Figure \ref{fig:conf-curve} for each of these 37 contrasts are available at \url{https://github.com/pilsneyrouset/comparison_Notip-pARI}.

\begin{figure}[!ht]
    \centering
    \includegraphics[width=12cm]{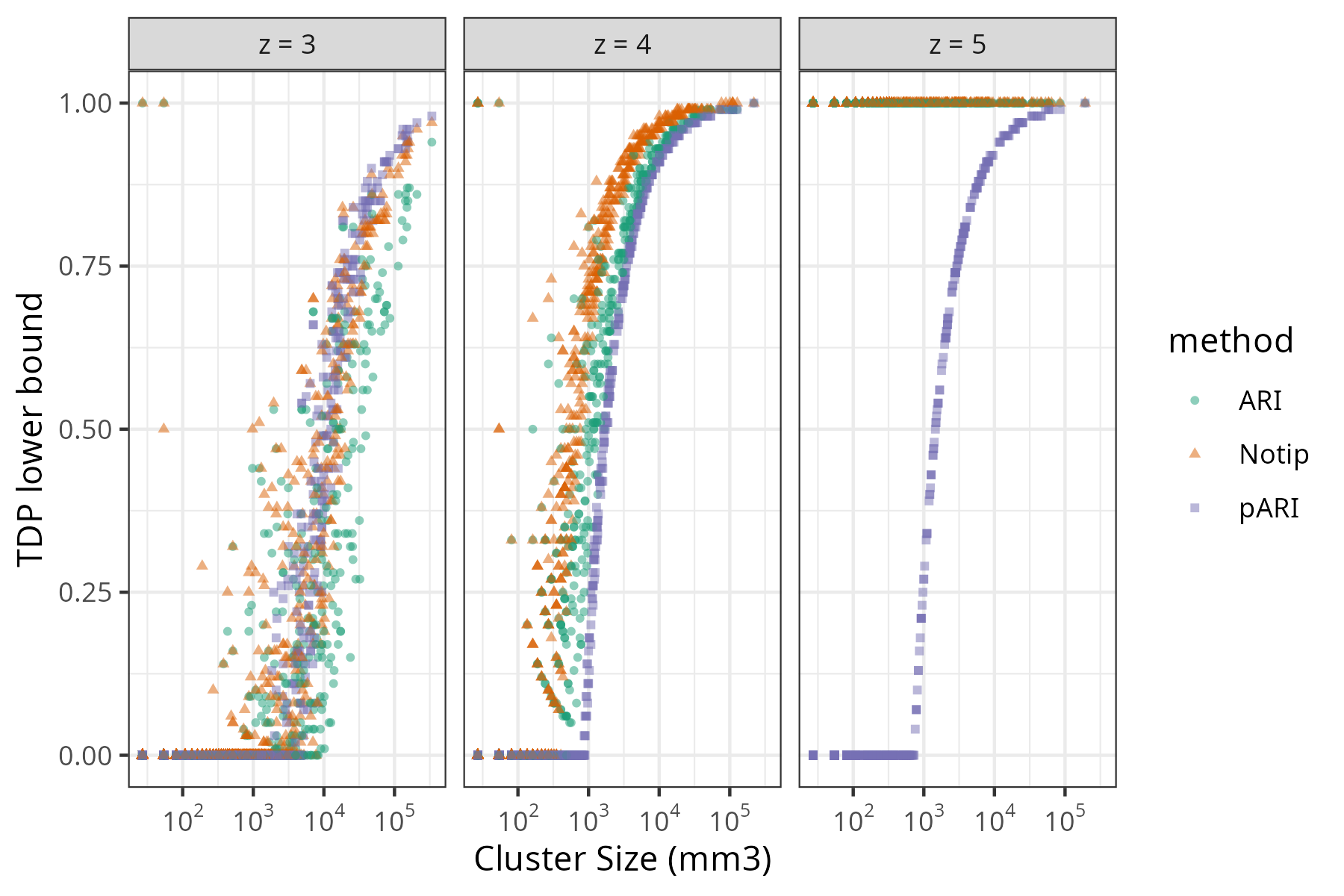}
    \caption{Lower bound on the True Discovery Proportion $\overline{\mathrm{TDP}}(S)$ as a function of cluster size $|S|$, for each cluster $S$ identified at the cluster forming threshold $z=3$ (left panel), $z=4$ (center panel), and $z=5$ (right panel).
}
    \label{fig:TDP-vs-cluster-size}
\end{figure}

These results are summarized in Figure~\ref{fig:TDP-vs-cluster-size}, where each panel corresponds to a value for the cluster-forming threshold $z$.
For each method, the TDP bound of each cluster is plotted against its size.
For a given value of $z$, larger clusters correspond to regions where the signal is stronger. Accordingly, the TDP bounds tend to be larger for larger clusters for a given method.

Comparing between methods based on the TDP bounds leads to the following conclusions. 
First, Notip consistently outperforms ARI. In particular, as noted above in the case of the dataset studied in Section~\ref{sec:cluster-eval-one-contrast}, it retains and improves upon the drill-down ability of the ARI method.

When the cluster-forming threshold is set very high ($z=5$), the signal is so strong that both methods report that the TDP is equal to 1. This  corresponds to pure signal in all $481$ clusters obtained across the 37 datasets: those are subsets of the Family-Wise Error Rate (FWER)-controlling set provided by the Bonferroni-Holm method \citep{holm1979simple}.
In contrast, the pARI method detects pure signal in only $1$ of these $481$ clusters.

For a small value of the cluster-forming threshold ($z=3$), pARI outperforms ARI for larger clusters, and its performance is globally similar to Notip's. 
Consistent with the results reported by \cite{andreella2024selective}, pARI slightly outperforms Notip for larger clusters. 
However, the results differ markedly for smaller clusters.  Here, pARI almost always yields null TDP bounds, while both Notip and ARI provide informative (i.e. non-null) TDP bounds.

Naturally, increasing the value of the cluster-forming threshold value to $z=4$ leads to larger TDP bounds for all methods. 
However this improvement is not uniform across methods. pARI  systematically underperforms compared to Notip and the baseline ARI method.
pARI's performance drops even more dramatically at $z=5$, where it provides uninformative TDP bounds for clusters of size below $1,000$ mm$^3$ and non-trivial but massively underestimated TDP bounds for larger clusters.

\section{Discussion}
\label{sec:discussion}
We have performed an extensive comparison between two recently proposed methods for post hoc inference for fMRI data: Notip \citep{blain2022notip} and pARI \citep{andreella2023permutation}. 
As expected from the theory, since both are based on the same calibration principle \citep{blanchard2020post}, our numerical experiments confirm that neither of Notip nor pARI is uniformly more powerful than the other (no free lunch).

This study illustrates the importance and difficulty of objectively assessing the performance of methods. 
The Notip paper \citep{blain2022notip} focused on the size of the largest detected regions because the relative behavior of the different compared methods did not depend on the the region size. Notip remains consistently better than the baseline method, ARI.
In contrast, the pARI method introduces a parameter $\delta$, which indirectly controls the size of the smallest cluster for which informative TDP bounds can be obtained \citep{andreella2023permutation}.
Therefore, performance comparisons involving pARI must also consider smaller clusters. 

Our experiments have shown that pARI can perform dramatically worse than Notip and even the baseline method ARI, especially in regions with a large amount of signal. 
Unfortunately, this precludes the drill-down approach advocated by \cite{rosenblatt2018ari}, where the cluster-forming threshold is increased in order to locate the signal with greater spatial precision.
This limitation can be problematic since low cluster forming thresholds have been shown to lead to unreliable inference \citep{Woo2014}. Specifically, \citep{Woo2014} ``recommend setting $p < .001$ as a lower limit default, and using more stringent primary thresholds or voxel-wise correction methods for highly powered studies''.

The methods discussed in this paper are not specific to fMRI studies and can be used in other contexts, such as genomics (see e.g. \cite{enjalbert-courrech2022powerful}), provided relevant hyperparameters are chosen. 

Finally, we would like to remind users that the hyperparameters discussed in this work ($\delta$ for pARI and $K$ for Notip) must be set prior to data analysis. Selecting them after the fact is another instance of double-dipping.

\bibliographystyle{abbrvnat}
\bibliography{posthoc}

\appendix
\section*{Appendix}

\section{Additional numerical results}
\subsection{``Look  negative  cue''  vs  ``Look  negative  rating''  dataset}
\label{sec:additional-tables}
In Table~\ref{tab:TDP-cluster-table_task36_z=3}, we reproduce the results obtained for $z=3$ in \cite[Table 2]{blain2022notip} and complemented by \cite[Table 1]{andreella2024selective}. 
Note that the results of Table \ref{tab:TDP-cluster-table_task36_z=3} are not exactly identical to those in \cite{andreella2024selective} because of the numerical variability inherent to the use of random permutation in the analyis.
\begin{table}[htp]
    \centering
    
\begin{tabular}{lrrrrrlll}
\toprule
\multicolumn{6}{c}{ } & \multicolumn{3}{c}{TDP lower bound} \\
\cmidrule(l{3pt}r{3pt}){7-9}
ID & X & Y & Z & Peak Stat & Size (mm3) & ARI & Notip & pARI\\
\midrule
1 & -33 & -94 & -17 & 5.63 & 7695 & 0.17 & 0.3 & \textbf{0.33}\\
2 & 66 & 2 & 16 & 5.47 & 14877 & 0.2 & 0.46 & \textbf{0.57}\\
3 & -12 & -82 & -8 & 5.40 & 14445 & 0.27 & 0.52 & \textbf{0.59}\\
4 & -6 & 11 & 52 & 5.30 & 5238 & 0.14 & 0.31 & \textbf{0.32}\\
5 & 45 & 14 & 25 & 5.27 & 4563 & 0.24 & \textbf{0.33} & 0.28\\
6 & 12 & -43 & -26 & 5.08 & 12555 & 0.05 & 0.36 & \textbf{0.5}\\
7 & 39 & -73 & 4 & 5.00 & 6075 & 0.04 & 0.21 & \textbf{0.23}\\
8 & -63 & -34 & 16 & 4.95 & 25812 & 0.3 & 0.66 & \textbf{0.75}\\
9 & 36 & -94 & -8 & 4.75 & 6507 & 0.08 & \textbf{0.19} & \textbf{0.19}\\
\bottomrule
\end{tabular}

    \caption{``Look  negative  cue''  vs  ``Look  negative  rating''  dataset: comparison between lower bounds on the True Discovery Proportion for the cluster-defining threshold $z=3$.}
    \label{tab:TDP-cluster-table_task36_z=3}
\end{table}
We also provide additional results corresponding to $z=4$, $4.5$ and $5$ in Tables~\ref{tab:TDP-cluster-table_task36_z=4}, Tables~\ref{tab:TDP-cluster-table_task36_z=4.5} and~\ref{tab:TDP-cluster-table_task36_z=5}. 

\begin{table}[htp]
    \centering
    
\begin{tabular}{lrrrrrlll}
\toprule
\multicolumn{6}{c}{ } & \multicolumn{3}{c}{TDP lower bound} \\
\cmidrule(l{3pt}r{3pt}){7-9}
ID & X & Y & Z & Peak Stat & Size (mm3) & ARI & Notip & pARI\\
\midrule
1 & -33 & -94 & -17 & 5.63 & 1431 & 0.64 & \textbf{0.74} & 0.42\\
2 & 66 & 2 & 16 & 5.47 & 2997 & 0.74 & \textbf{0.87} & 0.72\\
3 & -12 & -82 & -8 & 5.40 & 1431 & 0.58 & \textbf{0.75} & 0.42\\
4 & -6 & 11 & 52 & 5.30 & 1485 & 0.51 & \textbf{0.76} & 0.44\\
5 & 45 & 14 & 25 & 5.27 & 1755 & 0.62 & \textbf{0.78} & 0.52\\
6 & 12 & -43 & -26 & 5.08 & 459 & 0.35 & \textbf{0.47} & 0\\
7 & 39 & -73 & 4 & 5.00 & 405 & 0.2 & \textbf{0.4} & 0\\
8 & 30 & -73 & -8 & 4.96 & 567 & 0.29 & \textbf{0.43} & 0\\
9 & -63 & -34 & 16 & 4.95 & 3726 & 0.79 & \textbf{0.9} & 0.78\\
10 & -24 & -61 & -11 & 4.91 & 594 & 0.32 & \textbf{0.45} & 0\\
11 & -27 & -19 & 4 & 4.85 & 216 & \textbf{0.25} & \textbf{0.25} & 0\\
12 & 36 & -94 & -8 & 4.75 & 1134 & 0.48 & \textbf{0.69} & 0.26\\
13 & 30 & -46 & -11 & 4.64 & 1188 & 0.5 & \textbf{0.68} & 0.3\\
14 & -60 & -49 & 25 & 4.59 & 324 & 0 & \textbf{0.08} & 0\\
15 & -45 & -79 & -26 & 4.56 & 513 & 0 & \textbf{0.32} & 0\\
20 & 0 & -64 & -14 & 4.43 & 378 & 0 & \textbf{0.07} & 0\\
23 & -45 & -67 & 34 & 4.32 & 405 & 0 & \textbf{0.13} & 0\\
\bottomrule
\end{tabular}

    \caption{``Look  negative  cue''  vs  ``Look  negative  rating''  dataset: comparison between lower bounds on the True Discovery Proportion for the cluster-defining threshold $z=4$.}
    \label{tab:TDP-cluster-table_task36_z=4}
\end{table}

\begin{table}[htp]
    \centering
    
\begin{tabular}{lrrrrrlll}
\toprule
\multicolumn{6}{c}{ } & \multicolumn{3}{c}{TDP lower bound} \\
\cmidrule(l{3pt}r{3pt}){7-9}
ID & X & Y & Z & Peak Stat & Size (mm3) & ARI & Notip & pARI\\
\midrule
1 & -33 & -94 & -17 & 5.63 & 918 & 0.91 & \textbf{0.94} & 0.21\\
2 & 66 & 2 & 16 & 5.47 & 513 & \textbf{0.89} & \textbf{0.89} & 0\\
3 & -12 & -82 & -8 & 5.40 & 783 & \textbf{0.93} & \textbf{0.93} & 0.07\\
4 & -6 & 11 & 52 & 5.30 & 594 & 0.86 & \textbf{0.91} & 0\\
5 & 45 & 14 & 25 & 5.27 & 783 & \textbf{0.93} & \textbf{0.93} & 0.07\\
6 & 12 & -43 & -26 & 5.08 & 189 & \textbf{0.86} & \textbf{0.86} & 0\\
7 & 39 & -73 & 4 & 5.00 & 81 & \textbf{0.67} & \textbf{0.67} & 0\\
8 & 30 & -73 & -8 & 4.96 & 189 & \textbf{0.71} & \textbf{0.71} & 0\\
9 & -63 & -34 & 16 & 4.95 & 702 & 0.88 & \textbf{0.92} & 0\\
10 & -24 & -61 & -11 & 4.91 & 162 & \textbf{0.67} & \textbf{0.67} & 0\\
11 & -63 & -10 & 13 & 4.90 & 108 & \textbf{0.5} & \textbf{0.5} & 0\\
12 & -27 & -19 & 4 & 4.85 & 81 & \textbf{0.67} & \textbf{0.67} & 0\\
13 & 36 & -94 & -8 & 4.75 & 432 & 0.81 & \textbf{0.88} & 0\\
14 & -57 & -19 & 7 & 4.68 & 108 & \textbf{0.5} & \textbf{0.5} & 0\\
15 & 69 & -22 & 10 & 4.67 & 108 & \textbf{0.5} & \textbf{0.5} & 0\\
16 & 30 & -46 & -11 & 4.64 & 270 & 0.7 & \textbf{0.8} & 0\\
\bottomrule
\end{tabular}

    \caption{``Look  negative  cue''  vs  ``Look  negative  rating''  dataset: comparison between lower bounds on the True Discovery Proportion for the cluster-defining threshold $z=4.5$.}
    \label{tab:TDP-cluster-table_task36_z=4.5}
\end{table}

\begin{table}[htp]
    \centering
    
\begin{tabular}{lrrrrrlll}
\toprule
\multicolumn{6}{c}{ } & \multicolumn{3}{c}{TDP lower bound} \\
\cmidrule(l{3pt}r{3pt}){7-9}
ID & X & Y & Z & Peak Stat & Size (mm3) & ARI & Notip & pARI\\
\midrule
1 & -33 & -94 & -17 & 5.63 & 378 & \textbf{1} & \textbf{1} & 0\\
2 & 66 & 2 & 16 & 5.47 & 135 & \textbf{1} & \textbf{1} & 0\\
3 & -12 & -82 & -8 & 5.40 & 135 & \textbf{1} & \textbf{1} & 0\\
4 & -6 & 11 & 52 & 5.30 & 81 & \textbf{1} & \textbf{1} & 0\\
5 & 45 & 14 & 25 & 5.27 & 216 & \textbf{1} & \textbf{1} & 0\\
6 & 12 & -43 & -26 & 5.08 & 27 & \textbf{1} & \textbf{1} & 0\\
7 & 39 & -73 & 4 & 5.00 & 27 & \textbf{1} & \textbf{1} & 0\\
\bottomrule
\end{tabular}

    \caption{``Look  negative  cue''  vs  ``Look  negative  rating''  dataset: comparison between lower bounds on the True Discovery Proportion for the cluster-defining threshold $z=5$.}
    \label{tab:TDP-cluster-table_task36_z=5}
\end{table}

\end{document}